\documentclass[seceqn,doublespacing,final,12pt]{elsart1p}

\newcommand{\ind}{1\hspace{-.27em}\mbox{\rm l}}

\usepackage{natbib}

\usepackage{graphicx, amsopn, amssymb, amsmath, amsfonts,vmargin}
\journal{Statistics \& Probability Letters}
\begin{document}

\begin{frontmatter}



\title{Some Darling-Siegert relationships connected with random flights}

 \author[lable1]{V. Cammarota},
 \ead{valentina.cammarota@uniroma1.it}
 \author[lable2]{A. Lachal\corauthref{cor1}},
 \ead{aime.lachal@insa-lyon.fr}
 \corauth[cor1]{Corresponding author.}
\author[lable1]{E. Orsingher}
 \ead{enzo.orsingher@uniroma1.it}

 \address[lable1]{Dipartimento di Statistica, Probabilit\`a e Statistiche Applicate, University of Rome `La Sapienza', P.le Aldo Moro 5, 00185 Rome, Italy}
  \address[lable2]{P\^{o}le de Math\'ematiques, Institut National des Sciences Appliqu\'ees de Lyon, B\^atiment L\'eonard de Vinci, 20 avenue Albert Einstein, 69621 Villeurbanne Cedex, France}

\author{}

\address{}

\begin{abstract}
% Text of abstract
We derive in detail four important results on integrals of Bessel functions from which three combinatorial identities are extracted. We present the probabilistic interpretation of these identities in terms of different types of random walks, including asymmetric ones. This work extends the results of a previous paper concerning the Darling-Siegert interpretation of similar formulas emerging in the analysis of random flights.
\end{abstract}


\end{frontmatter}

% main text
\section{Introduction}
\label{}

The analysis of random flights (see \cite{O-D}) performed by a particle changing direction at Poisson paced times with uniformly distributed orientations of the steps lets emerge a number of intriguing relationships involving different types of Bessel integrals. These are the source of identities with curious combinations of Gamma functions which suggest interesting probabilistic interpretations.

The instrument to give a probabilistic background to these combinatorial formulas are random walks on the set $\mathbb{Z}$ (with upward and downward steps) or on $\mathbb{N}$ (with non-decreasing sample paths). We arrive at different forms of Darling-Siegert expressions where the distribution of first-passage times and transition functions are suitably combined. We work on generating functions (instead of the reflection principle as in \cite{D-S}) and obtain directly many distributions, including the first passage time through an arbitrary level $a$ when steps take values $\pm 1$ with probability $p$ and $q$. The case where some conditional first hitting time of random walks is related to combinatorial identities is also examined below.

In the analysis of a planar random motion with uniform orientation of the steps, the semigroup property of the Bessel functions (\ref{16:40A}) yields as a byproduct the identity (\ref{b}) which suggests an interpretation based on a non-decreasing random walk with steps taking values $0$ or $1$.

Our aim is to make the paper as self-contained as possible and for this reason all the Bessel integrals are proved in detail as well as the distributions used in the underlying probabilistic interpretation.

We recall the usual representations of Bessel functions (\cite{G-R}, formulas 8.402 and 8.411.1): for $\mu >0$, the function $J_\mu$ is defined by the following equivalent equalities
\begin{eqnarray}
J_\mu(x)&=&\sum_{n=0}^{\infty} \frac{(-1)^n}{n! \, \Gamma(n+\mu+1)} \left(\frac{x}{2}\right)^{\!2n+\mu}\!, \label{ser} \\
J_\mu(x)&=& \frac{1}{2 \pi} \int_{-\pi}^{\pi} e^{i x \sin \theta -i \mu \theta} \,\mathrm{d} \theta. \label{int}
\end{eqnarray}
We denote by $\mathbb{N}$ the set of the positive integer numbers including zero and by $\mathbb{Z}$ that of the integer numbers.

Our analysis starts with the following integrals of Bessel functions.
\begin{prop} \label{12:25}
For any $\mu, \nu, a >0$, we have that
\begin{eqnarray}
\int_{0}^{a} \frac{J_\mu(x)J_\nu(a-x)}{x} \,\mathrm{d}x&=&\frac{1}{\mu} \, J_{\mu+\nu}(a), \label{16:38} \\
\int_{0}^{a} \frac{J_\mu(x)J_\nu(a-x)}{x(a-x)} \,\mathrm{d}x&=&\left(\frac{1}{\mu}+ \frac{1}{\nu}\right) \frac{J_{\mu+\nu}(a)}{a}, \label{16:39} \\
\int_{0}^{a} x^\mu (a-x)^{\nu} J_\mu (x) J_\nu (a-x) \,\mathrm{d}x &=& \frac{\Gamma(\mu+1/2) \, \Gamma(\nu+1/2)}{\sqrt{2 \pi}\, \Gamma(\mu+\nu+1)} \, a^{\mu+\nu+1/2} J_{\mu+\nu+1/2} (a), \label{16:40A} \\
\int_{0}^{a} x^\mu (a-x)^{\nu} J_{\mu-1} (x) J_\nu (a-x) \,\mathrm{d}x &=& \frac{\Gamma(\mu+1/2) \, \Gamma(\nu+1/2)}{\sqrt{2 \pi}\, \Gamma(\mu+\nu+1)} \, a^{\mu+\nu+1/2} J_{\mu+\nu-1/2} (a). \label{16:40B}
\end{eqnarray}
\end{prop} \label{jhs}

Actually, the above formulas hold for larger domains of the parameters $\mu$ and $\nu$ (see \cite{G-R}), but we shall only deal with positive $\mu$ and $\nu$.

Let us emphasize that (\ref{16:38}), (\ref{16:39}), (\ref{16:40A}) and (\ref{16:40B}) are particular cases of the general family of Bessel integrals $\int_{0}^{a} x^\alpha (a-x)^\beta J_\mu(x) J_\nu(a-x) \,\mathrm{d}x$ and simple formulas exist only in very few cases.

From Proposition \ref{12:25} we deduce the following formulas for the Gamma function.
\begin{thm} \label{thm2}
For any $\mu, \nu >0$ and $r \in \mathbb{N}$, we have that
\begin{eqnarray}
\sum_{k=0}^{r} \frac{\Gamma(2k+\mu)}{k! \, \Gamma(k+\mu+1)} \, \frac{\Gamma(2r-2k+\nu)}{(r-k)! \,\Gamma(r-k+\nu)}&=&\frac{1}{\mu} \,\frac{\Gamma(2r+\mu+\nu)}{r! \, \Gamma(r+\mu+\nu)}, \label{s} \\
\sum_{k=0}^{r} \frac{\Gamma(2k+\mu)}{k! \, \Gamma(k+\mu+1)} \, \frac{\Gamma(2r-2k+\nu)}{(r-k)! \,\Gamma(r-k+\nu+1)}&=&\left(\frac{1}{\mu}+\frac{1}{\nu}\right)\frac{\Gamma(2r+\mu+\nu)}{r! \, \Gamma(r+\mu+\nu+1)}, \label{d} \\
\sum_{k=0}^{r} \frac{\Gamma(k+\mu)}{k!} \frac{\Gamma(r-k+\nu)}{(r-k)!}&=& \frac{\Gamma(\mu) \, \Gamma(\nu)}{\Gamma(\mu+\nu)} \,\frac{\Gamma(r+\mu+\nu)}{r!}. \label{b}
\end{eqnarray}
\end{thm}
Formulas (\ref{s}) and (\ref{d}) come respectively from (\ref{16:38}) and (\ref{16:39}) while (\ref{b}) comes from (\ref{16:40A}) or, equivalently, (\ref{16:40B}).

Formulas (\ref{s}) and (\ref{d}) lead to the following Darling-Siegert-type relationships (see \cite{D-Si} and also \cite{D-S} for the discrete-time version) involving the distribution of the classical random walk $\left(S_n\right)_{n\in \mathbb{N}}$ on $\mathbb{Z}$ with steps taking values $+1$ or $-1$. Let us introduce $T_a=\min\{n \ge 1: S_n=a\}$ for any $a \in \mathbb{Z}$, with the usual convention that $\min (\emptyset) = + \infty$. The random variable $T_a$ is the first hitting (discrete) time of level $a$ by the random walk.
\begin{cor} \label{cor 1}
For any $\mu, \nu, r \in \mathbb{N} \setminus \{0\}$, we have that
\begin{eqnarray}
\sum_{k=0}^{r} \mathbb{P}\{T_\mu=2k+\mu\} \, \mathbb{P}\{S_{2r-2k+\nu}=\nu\} &=&\mathbb{P}\{S_{2r+\mu+\nu}=\mu+\nu\}, \label{16:37}\\
\sum_{k=0}^{r} \mathbb{P}\{T_\mu=2k+\mu\} \, \mathbb{P}\{T_\nu=2r-2k+\nu\} &=&\mathbb{P}\{T_{\mu+\nu}=2r+\mu+\nu\}. \label{16:37'}
\end{eqnarray}
\end{cor}
Formulas (\ref{s}) and (\ref{d}), and then (\ref{16:37}) and (\ref{16:37'}), have already been stated in \cite{D-S}. However, therein, they were related to the symmetric random walk instead. It is remarkable that all random walks with steps $\pm 1$ lead to the universal identities (\ref{s}) and (\ref{d}) whatever the probability of jumps.

Similarly, let $\left(S'_n\right)_{n \in \mathbb{N}}$ be the random walk with steps taking values $0$ or $1$ and let $T'_a=\min\{n \ge 1: S'_n=a\}$. Formula (\ref{b}) leads to both following Darling-Siegert-type relationships.
\begin{cor}
For any $\mu, \nu, r \in \mathbb{N}$, we have that
\begin{eqnarray}
\sum_{k=0}^{r} \mathbb{P}\{T'_\mu=k+\mu\} \, \mathbb{P}\{S'_{r-k+\nu}=\nu\} &=&\mathbb{P}\{S'_{r+\mu+\nu}=\mu+\nu\}, \label{16:37A'}\\
\sum_{k=0}^{r} \mathbb{P}\{T'_\mu=k+\mu\} \, \mathbb{P}\{T'_\nu=r-k+\nu\} &=&\mathbb{P}\{T'_{\mu+\nu}=r+\mu+\nu\}. \label{16:37A}
\end{eqnarray}
\end{cor}
It is noteworthy that (\ref{16:37A'}) and (\ref{16:37A}) come from the only universal identity (\ref{b}). This fact hinges on a feature of the non-decreasing random walk $\left(S'_n\right)_{n \in \mathbb{N}}$ and will be explained in Subsection \ref{subsection} and Remark \ref{remark}. For the same reason, (\ref{16:40A}) and (\ref{16:40B}) are equivalent.

Of course, formulas (\ref{16:37}), (\ref{16:37'}), (\ref{16:37A'}) and (\ref{16:37A}) can be derived by means of independent arguments by using the Markov property of random walks as explained in Subsection \ref{back}. We have chosen to prove them by means of Theorem \ref{thm2} (or equivalently of Proposition \ref{12:25}) because of the importance of Bessel funtions in many fields of Mathematics and especially in Probability Theory. There is a very huge amount of formulas concerning the special functions of Mathematical Physics, and we find it interesting to provide for certain ones, when possible, a probabilistic interpretation.

Our paper is organized as follows. In Section \ref{2}, we recall some well-known formulas of Bessel integrals. Some proofs of them relying on hypergeometrical functions can be found in \cite{watson}, chapters 12 and 13. Nevertheless, elementary proofs are not so accessible in the literature and then, because of their interest, we provide them in order to make the paper self-contained and easily readable. In Section \ref{tree}, we derive the Darling-Siegert-type interpretation of formulas (\ref{s}), (\ref{d}) and (\ref{b}). We point out that the Darling-Siegert-type relationship appearing in \cite{D-Si} deals with continuous Markov processes whereas our work concerns discrete Markov processes.

The main tool of our analysis is the use of generating functions for the probabilities appearing in Corollary \ref{cor 1}.

\section{Mathematical Background} \label{2}

Let us recall well-known formulas which will be used throughout the paper:
\begin{equation} \label{duplic}
\mbox{for $x>0$,} \quad \Gamma(2x)=\frac{2^{2x-1}}{\sqrt{\pi}} \, \Gamma(x) \, \Gamma(x+1/2),
\end{equation}
\begin{equation} \label{serie}
\mbox{for $|x|<1$,} \quad (1-x)^{-a}=\left\{\begin{array}{ll} \displaystyle \sum_{n=0}^{\infty} \frac{\Gamma(a+n)}{\Gamma(a)\,n!}\,x^n & \mbox{for $a>0$,} \\ \displaystyle \sum_{n=0}^{\infty} \binom{n+a-1}{a-1} x^n & \mbox{for $a\in \mathbb{N}\setminus\{0\}$,} \end{array}\right.
\end{equation}
\begin{equation} \label{seriepart}
\mbox{for $|x|<1$,} \quad (1-x)^{-1/2}= \sum_{n=0}^{\infty} \binom{2n}{n} \left(\frac{x}{4}\right)^{\!n}, \qquad (1-x)^{1/2}=1-2 \sum_{n=1}^{\infty} \frac{1}{n} \binom{2n-2}{n-1} \left(\frac{x}{4}\right)^{\!n}.
\end{equation}

%%%%%%%%%%%%%%%%%%%%%%%%%%%%%
\subsection{Proof of Proposition \ref{12:25}}
The first step for deriving formulas (\ref{16:38}) and (\ref{16:39}) is the evaluation of the Laplace transforms of the Bessel function $J_\nu$ which is presented in the next lemma (see \cite{G-R} formulas 6.611.1 p. 707 and 6.621.1 p. 711).
\begin{lem} \label{14:13}
For $\alpha, \beta, \nu >0$, we have that
\begin{eqnarray}
&&\int_{0}^{\infty}e^{-\alpha x} J_\nu (\beta x) \,\mathrm{d}x=\frac{(\sqrt{\alpha^2+\beta^2}-\alpha)^\nu}{\beta^\nu \sqrt{\alpha^2+\beta^2}}, \label{19:45cioccolata} \\
&&\int_{0}^{\infty}e^{-\alpha x} \frac{J_\nu (\beta x)}{x} \,\mathrm{d}x= \frac{(\sqrt{\alpha^2+\beta^2}-\alpha)^\nu}{\nu \beta^\nu}, \label{19:46americani} \\
&&\int_{0}^{\infty}e^{-\alpha x} x^\nu J_\nu(\beta x) \,\mathrm{d}x=\frac{2^\nu \Gamma(\nu+1/2) \, \beta^\nu}{\sqrt{\pi} \, (\alpha^2+\beta^2)^{\nu+1/2}}, \label{19:47}\\
&&\int_{0}^{\infty}e^{-\alpha x} x^\nu J_{\nu-1}(\beta x) \,\mathrm{d}x=\frac{2^\nu \Gamma(\nu+1/2) \, \alpha\beta^{\nu-1}}{\sqrt{\pi}\, (\alpha^2+\beta^2)^{\nu+1/2}}. \label{19:47A}
\end{eqnarray}
\end{lem}
For integer values of the index $\nu$ we are able to present a proof of these formulas.

%%%%%%%%%%%%%%%%%%%%%%%%%%%%%
\begin{pf}\\
%%%%%%%%%%%%%%%%%%%%%%%%%%%%%
\noindent $\bullet$ \textsl{Proof of (\ref{19:45cioccolata})}. We consider the integral representation (\ref{int}) of the Bessel function which we rewrite as
$
J_\nu(x)=\frac{1}{2 \pi} \int_{-\pi}^{\pi} e^{i \beta x \sin \theta -i \nu \theta} \,\mathrm{d} \theta
= \frac{e^{-i \nu \frac{\pi}{2}}}{2 \pi} \int_{0}^{2 \pi} e^{i \beta x \cos \theta + i \nu \theta} \,\mathrm{d} \theta
$
and we get
$$
\int_{0}^{\infty} e^{-\alpha x} J_\nu(\beta x) \,\mathrm{d}x = \frac{e^{-i \nu \frac{\pi}{2}}}{2 \pi} \int_{0}^{2\pi} e^{i \nu \theta} \,\mathrm{d} \theta \int_{0}^{\infty} e^{-(\alpha-i \beta \cos \theta)x} \,\mathrm{d}x = \frac{e^{-i \nu \frac{\pi}{2}}}{2 \pi} \int_{0}^{2 \pi} \frac{e^{i \nu \theta}}{\alpha-i \beta \cos \theta} \,\mathrm{d} \theta.
$$
For any integer number $\nu$, the function $f(z)= \frac{e^{i \nu z}}{\alpha-i \beta \cos z}$ is periodic with period $2 \pi$ and its poles are $\pm \left(\frac{\pi}{2} + i \arg\sinh \frac{\alpha} {\beta}\right) +2 k \pi$, $k \in \mathbb{Z}$. We take a contour $\Gamma_R$ consisting of $\Gamma_1=\{(x,y): x \in [0, 2 \pi], y=0\}$, $\Gamma_2=\{(x,y): x=2 \pi, y \in [0,R]\}$, $\Gamma_3=\{(x,y): x \in [0, 2 \pi], y=R\}$ and $\Gamma_4=\{(x,y): x=0, y \in [0,R]\}$ with $R$ such that $\Gamma_R$ encloses the unique pole $\theta_0=\frac{\pi}{2}+ i \arg \sinh \frac{\alpha}{\beta}$ with residue $\rho_0=\frac{e^{i \nu \theta_0}}{i \beta \sin \theta_0}=e^{i \nu \frac{\pi}{2}} \frac{( \sqrt{\alpha^2+\beta^2}-\alpha )^\nu}{i \beta^\nu \sqrt{\alpha^2+\beta^2} }$. By Cauchy's residue theorem we obtain (\ref{19:45cioccolata}).

%%%%%%%%%%%%%%%%%%%%%%%%%%%%%
\noindent $\bullet$ \textsl{Proof of (\ref{19:46americani})}. It is straightforward to derive (\ref{19:46americani}) from (\ref{19:45cioccolata}) as follows:
\begin{eqnarray}
\int_{0}^{\infty} e^{-\alpha x} \frac{J_\nu(\beta x)}{x} \,\mathrm{d}x=\int_{\alpha}^{\infty} \left(\int_{0}^{\infty} e^{-a x} J_\nu(\beta x) \,\mathrm{d}x\right) \mathrm{d} a=\frac{(\sqrt{\alpha^2+\beta^2}-\alpha)^\nu}{\nu \beta^\nu}. \nonumber
\end{eqnarray}

%%%%%%%%%%%%%%%%%%%%%%%%%%%%%
\noindent $\bullet$ \textsl{Proof of (\ref{19:47})}. In order to prove (\ref{19:47}), we use the equality
\begin{equation} \label{19:04}
\frac{\mathrm{d}}{\mathrm{d}x} \left(x^{\nu} J_{\nu}(\beta x)\right)=\beta x^{\nu} J_{\nu-1} (\beta x)
\end{equation}
which can be obtained by a direct calculation. Set $U_\nu=\int_{0}^{\infty} e^{-\alpha x} x^\nu J_\nu(\beta x) \,\mathrm{d}x$. In light of (\ref{19:04}), after an integration by parts, we obtain that
\begin{equation} \label{19:58}
U_\nu = \frac{1}{\alpha} \int_{0}^{\infty} e^{-\alpha x} \frac{\mathrm{d}}{\mathrm{d}x} \left(x^{\nu} J_{\nu}(\beta x)\right) \mathrm{d}x = \frac{\beta}{\alpha} \int_{0}^{\infty} e^{-\alpha x} x^{\nu} J_{\nu-1}(\beta x) \,\mathrm{d}x.
\end{equation}
A further integration by parts and formulas (\ref{19:04}) and (\ref{19:58}) together with $xe^{-\alpha x}=-\frac{\mathrm{d}}{\mathrm{d}x}((\frac{x}{\alpha}+\frac{1}{\alpha^2})e^{-\alpha x})$ give
\begin{eqnarray} \label{19:498}
U_\nu &=&\frac{\beta}{\alpha} \int_{0}^{\infty} (xe^{-\alpha x}) (x^{\nu-1} J_{\nu-1}(\beta x)) \,\mathrm{d}x =\frac{\beta}{\alpha} \int_{0}^{\infty} e^{-\alpha x} \left(\frac{x}{\alpha}+\frac{1}{\alpha^2}\right) \frac{\mathrm{d}}{\mathrm{d}x}( x^{\nu-1} J_{\nu-1}(\beta x) ) \,\mathrm{d}x \nonumber \\
&=&\frac{\beta^2}{\alpha^3} \int_{0}^{\infty} e^{-\alpha x} x^{\nu-1} J_{\nu-2} (\beta x) \,\mathrm{d}x+ \frac{\beta^2}{\alpha^2} \int_{0}^{\infty} e^{-\alpha x} x^{\nu} J_{\nu-2} (\beta x) \,\mathrm{d}x \nonumber \\
&=&\frac{\beta}{\alpha^2} \, U_{\nu-1}+\left(\frac{2(\nu-1)\beta}{\alpha^2}\, U_{\nu-1} -\frac{\beta^2}{\alpha^2}\, U_{\nu}\right) =\frac{(2\nu-1)\beta}{\alpha^2}\, U_{\nu-1} -\frac{\beta^2}{\alpha^2}\, U_{\nu}
\end{eqnarray}
where we made use of the relationship $J_{\nu}(\beta x)=\frac{2 (\nu-1)}{\beta x} J_{\nu-1}(\beta x)-J_{\nu-2}(\beta x)$ which can be proved by observing that $\frac{2 (\nu-1)}{\beta x} J_{\nu-1}(\beta x)-J_{\nu-2}(\beta x)=\sum_{n=0}^{\infty} \frac{(-1)^n}{n! \, \Gamma(n+\nu-1)} \left(\frac{\nu-1}{n+\nu-1} -1\right) \left(\frac{\beta x}{2}\right)^{\!2n+\nu-2}=J_\nu (\beta x) $. From (\ref{19:498}),  (\ref{duplic}) and $U_0=\frac{1}{\sqrt{\alpha^2+\beta^2}}$, we extract the recurrence relationship
$$
U_{\nu}=\frac{(2 \nu-1)\beta}{\alpha^2+\beta^2} \,U_{\nu-1}=\ldots=\frac{(2 \nu-1)(2\nu-3)\cdots 3 \, \beta^\nu}{(\alpha^2+\beta^2)^{\nu+1/2}}=\frac{2^{\nu} \Gamma(\nu+1/2) \, \beta^\nu}{\sqrt{\pi}(\alpha^2+\beta^2)^{\nu+1/2}}.
$$

%%%%%%%%%%%%%%%%%%%%%%%%%%%%%
\noindent $\bullet$ \textsl{Proof of (\ref{19:47A})}. From (\ref{19:58}), we derive $\int_{0}^{\infty} e^{-\alpha x} x^\nu J_{\nu-1}(\beta x) \,\mathrm{d}x = \frac{\alpha}{\beta} \int_{0}^{\infty} e^{-\alpha x} x^{\nu} J_{\nu}(\beta x) \,\mathrm{d}x$ which, in light of (\ref{19:47}), yields (\ref{19:47A}).
\end{pf}

We now employ the Laplace transforms (\ref{19:45cioccolata}), (\ref{19:46americani}), (\ref{19:47}) and (\ref{19:47A}) to prove Proposition \ref{12:25} which plays an important role in our paper.

%%%%%%%%%%%%%%%%%%%%%%%%%%%%%
\vspace{\baselineskip} \noindent {\bf{\uppercase{Proof of proposition \ref{12:25}}}}.\\
%\begin{pf}\\
%%%%%%%%%%%%%%%%%%%%%%%%%%%%%
\noindent $\bullet$ \textsl{Proof of (\ref{16:38}) and (\ref{16:39})}. We can prove formulas (\ref{16:38}) and (\ref{16:39}) by taking the Laplace transforms of both members or directly as shown in \cite{D-S}. Let us evaluate the Laplace transforms of the left-hand sides of (\ref{16:38}) and (\ref{16:39}). Remarking that we are dealing with convolution products, we easily get, by Lemma \ref{14:13},
\begin{eqnarray*}
\int_{0}^{\infty} e^{- \alpha a}\,\mathrm{d}a \int_{0}^{a} \frac{J_\mu(x)J_\nu(a-x)}{x} \,\mathrm{d}x&=& \frac{(\sqrt{\alpha^2+\beta^2}-\alpha)^{\mu+\nu}}{\mu \beta^{\mu+\nu} \sqrt{\alpha^2+\beta^2}} =\frac{1}{\mu} \int_{0}^{\infty} e^{- \alpha a} J_{\mu+\nu}(a) \,\mathrm{d}a,
\end{eqnarray*}
and
\begin{eqnarray*}
\int_{0}^{\infty} e^{- \alpha a} \,\mathrm{d}a \int_{0}^{a} \frac{J_\mu(x)J_\nu(a-x)}{x(a-x)} \,\mathrm{d}x&=&\frac{(\sqrt{\alpha^2+\beta^2}-\alpha)^{\mu+\nu}}{\mu \nu \beta^{\mu+\nu}}=\frac{\mu+\nu}{\mu \nu} \int_{0}^{\infty} e^{-\alpha a} \frac{J_{\mu+\nu(a)}}{a} \,\mathrm{d}a.
\end{eqnarray*}
As a result, formulas (\ref{16:38}) and (\ref{16:39}) immediately emerge by inverting the Laplace transforms.

%%%%%%%%%%%%%%%%%%%%%%%%%%%%%
\noindent $\bullet$ \textsl{Proof of (\ref{16:40A}) and (\ref{16:40B})}. We can prove formulas (\ref{16:40A}) and (\ref{16:40B}) in the same way. Indeed, e.g., for (\ref{16:40A}),
\begin{eqnarray*}
\lefteqn{\int_{0}^{\infty} e^{- \alpha a} \,\mathrm{d}a \int_{0}^{a} x^\mu (a-x)^{\nu} J_\mu(x) J_\nu(a-x) \,\mathrm{d}x} \\
&=& \frac{2^{\mu+\nu} \Gamma(\mu+1/2) \, \Gamma(\nu+1/2)}{\pi (\alpha^2+1)^{\mu+\nu+1}} = \frac{\Gamma(\mu+1/2) \, \Gamma(\nu+1/2)}{\sqrt{2 \pi} \, \Gamma(\mu+\nu+1)} \int_{0}^{\infty} e^{- \alpha a} a^{\mu+\nu+1/2} J_{\mu+\nu+1/2}(a) \,\mathrm{d}a
\end{eqnarray*}
and thus, by inverting the Laplace transform, we conclude that (\ref{16:40A}) holds. The proof of (\ref{16:40B}) is quite similar.

\subsection{Proof of Theorem \ref{thm2}} \label{proof}

We now proceed to prove Theorem \ref{thm2}.

%%%%%%%%%%%%%%%%%%%%%%%%%%%%%%%%%%%%%%%%%%%%%%%%%%%%%%%%%%%%%%%%%%%%%%%%%%%%%%%%
\noindent $\bullet$ \textsl{Proof of (\ref{s})}. We start from (\ref{16:38}), using the series representation (\ref{ser}) and solving the integral $ \int_{0}^{a} x^{2k+\mu-1} (a-x)^{2l+\nu} \,\mathrm{d}x$, we can write that
\begin{eqnarray}
\lefteqn{\int_{0}^{a}\frac{J_\mu(x) J_\nu(a-x)}{x} \,\mathrm{d}x} \nonumber \\
&=& \sum_{k=0}^{\infty} \sum_{l=0}^{\infty} \frac{(-1)^{k+l}}{k! \, l!} \, \frac{\Gamma(2k+\mu) \, \Gamma(2l+\nu+1)}{\Gamma(k+\mu+1) \, \Gamma(l+\nu+1) \, \Gamma(2(k+l)+\mu+\nu+1)} \, \frac{a^{2(k+l)+\mu+\nu}}{2^{2(k+l)+\mu+\nu}} \nonumber \\
&=& \sum_{r=0}^{\infty} \frac{(-1)^r}{\Gamma(2r+\mu+\nu+1)} \left(\frac{a}{2}\right)^{\!2r+\mu+\nu} \sum_{k=0}^{r} \frac{\Gamma(2k+\mu)}{k! \, \Gamma(k+\mu+1)} \, \frac{\Gamma(2(r-k)+\nu+1)}{(r-k)! \, \Gamma(r-k+\nu+1)} \label{12:50}
\end{eqnarray}
and
\begin{equation} \label{12:51}
\frac{J_{\mu+\nu}(a)}{\mu}=\frac{1}{\mu} \sum_{r=0}^{\infty} \frac{(-1)^r}{r! \,\Gamma(r+\mu+\nu+1)} \left(\frac{a}{2}\right)^{\!2r+\mu+\nu}\!.
\end{equation}
By comparing the coefficients of the entire series (\ref{12:50}) and (\ref{12:51}), we immediately obtain formula (\ref{s}).

%%%%%%%%%%%%%%%%%%%%%%%%%%%%%%%%%%%%%%%%
\noindent $\bullet$ \textsl{Proof of (\ref{d})}. Analogously, by (\ref{16:39}), we write that
\begin{eqnarray}
\lefteqn{\int_{0}^{a}\frac{J_\mu(x)J_\nu(a-x)}{x(a-x)}\,\mathrm{d}x} \nonumber \\
&=&\sum_{k=0}^{\infty} \sum_{l=0}^{\infty} \frac{(-1)^{k+l}}{k!\, l!} \, \frac{\Gamma(2k+\mu) \, \Gamma(2l+\nu)}{\Gamma(k+\mu+1) \, \Gamma(l+\nu+1) \, \Gamma(2(k+l)+\mu+\nu)} \, \frac{a^{2(k+l)+\mu+\nu-1}}{2^{2(k+l)+\mu+\nu}} \nonumber \\
&=&\frac 1a\sum_{r=0}^{\infty} \frac{(-1)^r}{\Gamma(2r+\mu+\nu)} \left(\frac{a}{2}\right)^{\!2r+\mu+\nu} \sum_{k=0}^{r} \frac{\Gamma(2k+\mu) \, \Gamma(2r-2k+\nu)}{k! \, (r-k)!\,\Gamma(k+\mu+1) \, \Gamma(r-k+\nu+1)} \label{15:17}
\end{eqnarray}
and
\begin{equation} \label{15:17II}
\left(\frac{1}{\mu}+\frac{1}{\nu}\right) \frac{J_{\mu+\nu}(a)}{a}= \left(\frac{1}{\mu}+\frac{1}{\nu}\right) \frac{1}{a} \sum_{r=0}^{\infty} \frac{(-1)^r}{r! \,\Gamma(r+\mu+\nu+1)} \left(\frac{a}{2}\right)^{\!2r+\mu+\nu}\!.
\end{equation}
By comparing the coefficients of the entire series (\ref{15:17}) and (\ref{15:17II}), we obtain (\ref{d}).

%%%%%%%%%%%%%%%%%%%%%%%%%%%%%
\noindent $\bullet$ \textsl{Proof of (\ref{b})}. Similarly, from (\ref{16:40A}) and (\ref{duplic}), we have
\begin{eqnarray*}
\hspace{-0.4cm}\lefteqn{\int_{0}^{a} x^\mu (a-x)^\nu J_\mu(x) J_\nu(a-x) \,\mathrm{d}x} \\
&=&\sum_{k=0}^{\infty} \sum_{l=0}^{\infty} \frac{(-1)^{k+l}}{k!\, l!} \, \frac{\Gamma(2k+2\mu+1) \, \Gamma(2l+2\nu+1)}{\Gamma(k+\mu+1) \, \Gamma(l+\nu+1) \, \Gamma(2(k+l)+2(\mu+\nu)+2)} \, \frac{a^{2(k+l)+2(\mu+\nu)+1}}{2^{2(k+l)+\mu+\nu}} \\
&=&\frac{1}{\sqrt\pi} \sum_{k=0}^{\infty} \sum_{l=0}^{\infty} \frac{(-1)^{k+l}}{k!\, l!} \, \frac{\Gamma(k+\mu+1/2) \, \Gamma(l+\nu+1/2)}{\Gamma(k+l+\mu+\nu+1) \, \Gamma(k+l+\mu+\nu+3/2)} \, \frac{a^{2(k+l)+2(\mu+\nu)+1}}{2^{2(k+l)+\mu+\nu+1}} \\
&=&\frac{1}{\sqrt\pi} \sum_{r=0}^{\infty} \frac{(-1)^r}{\Gamma(r+\mu+\nu+3/2)} \frac{a^{2r+2(\mu+\nu)+1}}{2^{2r+\mu+\nu+1}} \sum_{k=0}^{r} \frac{\Gamma(k+\mu+1/2) \, \Gamma(r-k+\nu+1/2)}{k! \, (r-k)!\,\Gamma(r+\mu+\nu+1)}
\end{eqnarray*}
and
$$
\hspace{-0.1cm}\frac{\Gamma(\mu+1/2) \, \Gamma(\nu+1/2)}{\sqrt{2 \pi}\, \Gamma(\mu+\nu+1)} \,\frac{J_{\mu+\nu+1/2} (a)}{a^{-\mu-\nu-1/2}}
=\frac{1}{\sqrt{\pi}} \sum_{r=0}^{\infty} \frac{(-1)^r \, \Gamma(\mu+1/2) \, \Gamma(\nu+1/2)}{r! \, \Gamma(\mu+\nu+1) \, \Gamma(r+\mu+\nu+3/2)} \, \frac{a^{2r+2(\mu+\nu)+1}}{2^{2r+\mu+\nu+1}}.
$$
Then
$$
\sum_{k=0}^{r} \frac{\Gamma(k+\mu+1/2) \, \Gamma(r-k+\nu+1/2)}{k! \, (r-k)! \, \Gamma(r+\mu +\nu+1)} =\frac{\Gamma(\mu+1/2) \, \Gamma(\nu+1/2)}{r! \, \Gamma(\mu+\nu+1)}
$$
which proves, by translating $\mu$ and $\nu$ by $1/2$, formula (\ref{b}).

%%%%%%%%%%%%%%%%%%%%%%%%%%%%%
\section{Darling-Siegert-type relationship for random walks} \label{tree}

Formulas (\ref{s}), (\ref{d}) and (\ref{b}) can be interpreted by means of random walks as we explain below.

%%%%%%%%%%%%%%%%%%%%%%%%%%%%%
\subsection{Background on random walks} \label{back}

In order to facilitate the reading and because of its importance in many fields, we recall some well-known results on random walks on ${\mathbb Z}$ (see, e.g., \cite{feller} Chapter 3) together with their proofs. Fix $p\in (0,1)$, $q=1-p$ and let $\left(S_n\right)_{n \in \mathbb{N}}$ be a random walk on ${\mathbb Z}$ starting at the origin (i.e. $S_0=0$) with parameter $p$. For any $j\in \mathbb{Z}$, we introduce the generating function $G(\xi,j)$ of the probabilities $\mathbb{P}\{S_n=j\}$, $n \in \mathbb{N}$:
$$
G(\xi,j)=\sum_{n=0}^{\infty} \mathbb{P} \{S_n=j\} \xi^n.
$$
Remarking that the random variable $(S_n+n)/2$ has a binomial distribution with parameters $(n,p)$, we have, for $|j|\le n$ such that $n+j$ is even, that
\begin{equation} \label{17:42}
\mathbb{P}\{S_n=j\}= \binom{n}{\frac{n+j}{2}} p^{\frac{n+j}{2}} q^{\frac{n-j}{2}},
\end{equation}
and this probability vanishes when $n+j$ is odd. Notice the following equality which will be used further:
\begin{equation} \label{qqq}
\mathbb{P} \{S_n=-j\}= \left(\frac{q}{p}\right)^{\!j} \mathbb{P}\{S_n=j\}.
\end{equation}

\begin{lem}
For $|\xi|<1$ and $j \in \mathbb{Z}$, the function $G(\xi,j)$ is given by
\begin{equation} \label{21:41}
G(\xi,j)= \left\{\begin{array}{ll} \displaystyle{\frac{1}{\sqrt{1-4pq \xi^2}} \left(\frac{1- \sqrt{1-4pq \xi^2}}{2q \xi}\right)^{\!j}}& \mbox{for $j \ge 0$,} \\[3ex] \displaystyle{\frac{1}{\sqrt{1-4pq \xi^2}}\left(\frac{1- \sqrt{1-4pq \xi^2}}{2p \xi}\right)^{\!|j|}}& \mbox{for $j \le 0$.} \end{array}\right.
\end{equation}
\end{lem}
\begin{pf}
We can directly calculate $G(\xi,0)$ and $G(\xi,1)$: since $\mathbb{P}\{S_n=0\}$ vanishes for odd values of $n$, we have, by using (\ref{seriepart}), that
\begin{equation*}
G(\xi,0) = \sum_{n=0}^{\infty} \mathbb{P}\{S_n=0\} \xi^n=\sum_{n=0}^{\infty} \binom{2n}{n} (pq)^n \xi^{2n} = (1-4pq\xi^2)^{-1/2}.
\end{equation*}
Analogously, since $\mathbb{P}\{S_n=1\}$ vanishes for even values of $n$, we have, by using again (\ref{seriepart}), that
\begin{eqnarray*}
\hspace{-2.3em}G(\xi,1)&=& \sum_{n=1}^{\infty} \mathbb{P}\{S_n=1\} \xi^n = \sum_{n=0}^{\infty} \!\binom{2n+1}{n+1} p^{n+1} q^n \xi^{2n+1}=\sum_{n=0}^{\infty} \!\left[2 \binom{2n}{n}\!-\frac{1}{n+1} \binom{2n}{n}\!\right]\! p^{n+1} q^n \xi^{2n+1} \\
&=& \frac{2 p \xi}{\sqrt{1-4 p q \xi^2}} -\frac{1}{2 q \xi}\left(1-\sqrt{1-4 p q \xi^2}\,\right) = \frac{1-\sqrt{1-4pq\xi^2}}{2q\xi \sqrt{1-4pq\xi^2}}.
\end{eqnarray*}
Now, from the elementary recurrence relationship
$
G(\xi,j)= p \xi \, G(\xi,j-1)+ q \xi \, G(\xi, j+1),
$
which comes from the obvious identity $\mathbb{P}\{S_{n+1}=j\}=p \, \mathbb{P}\{S_n=j-1\}+q \, \mathbb{P}\{S_n=j+1\}$, we are able to prove (\ref{21:41}) by induction. First, for $j\ge 2$, by supposing that (\ref{21:41}) holds for $G(\xi,j-2)$ and $G(\xi,j-1)$, we have that
\begin{eqnarray*}
G(\xi,j)&=& \frac{1}{q \xi} \,G(\xi,j-1)- \frac{p}{q} \, G(\xi,j-2) \\
&=&\left(\frac{1-\sqrt{1-4pq\xi^2}}{2q\xi}\right)^{\!j-2} \frac{1}{\sqrt{1-4pq\xi^2}} \frac{1-4pq\xi^2-2 \sqrt{1-4pq\xi^2}+1}{4q^2\xi^2} \\
&=&\left(\frac{1-\sqrt{1-4pq \xi^2}}{2q\xi}\right)^{\!j} \frac{1}{\sqrt{1-4pq \xi^2}}.
\end{eqnarray*}
Therefore, (\ref{21:41}) holds for any positive integers $j$. Finally, we can easily check, by using (\ref{qqq}), that (\ref{21:41}) holds also for negative values of $j$: in fact, we can see that $G(\xi,-j)=(\frac{q}{p})^j G(\xi,j)$.
\end{pf}
Set now $T_a=\min\{n \ge 1: S_n=a\}$ for any integer $a$ with the usual convention $\min (\emptyset) =+\infty$. For $a \ne 0$, the random variable $T_a$ is the first hitting time of level $a$ for the random walk $(S_n)_{n \in \mathbb{N}}$ and $T_0$ is the fist return time to level $0$ for $(S_n)_{n \in \mathbb{N}}$.

\begin{prop} \label{Darling} (Darling-Siegert formulas) If $a$ and $b$ are both positive or both negative integer numbers, we have that
$$\mathbb{P}\{S_n=a+b\} = \sum_{k=a}^{n-b} \mathbb{P}\{T_a=k\} \, \mathbb{P}\{S_{n-k}=b\}, \hspace{0.5cm}\mathbb{P}\{T_{a+b}=n\}=\sum_{k=a}^{n-b}\mathbb{P}\{T_a=k\} \, \mathbb{P}\{T_b=n-k\}.$$
\end{prop}
\begin{pf}
Let $a \ge 1$ and $b \ge 0$ and let us consider a path connecting level $0$ at time $0$ to level $a+b$ at time $n$. This path must necessarily pass by level $a$ between times $0$ and $n$, that is, $T_a \le n$. So, from the strong Markov property of random walks, the translation invariance, and the fact that $S_{T_a}=a$ (in words, discrete continuity), we have
\begin{eqnarray*}
\mathbb{P}\{S_n=a+b\}&=& \sum_{k=0}^{n} \mathbb{P}\{T_a=k,\,S_n-a=b\} = \sum_{k=0}^{n} \mathbb{P}\{T_a=k\} \, \mathbb{P}\{S_n-S_k=b\} \\
&=& \sum_{k=0}^{n} \mathbb{P}\{T_a=k\} \, \mathbb{P}\{S_{n-k}=b\} = \sum_{k=a}^{n-b} \mathbb{P}\{T_a=k\} \, \mathbb{P}\{S_{n-k}=b\}
\end{eqnarray*}
and
\begin{eqnarray*}
\mathbb{P}\{T_{a+b}=n\}&=&\mathbb{P}\{T_a \le n, \, T_{a+b}=n\}=\sum_{k=0}^{n}\mathbb{P}\{T_a=k,\, T_{a+b}=n\} \\
&=&\sum_{k=0}^{n}\mathbb{P}\{T_a=k\} \, \mathbb{P}\{T_b=n-k\} =\sum_{k=a}^{n-b}\mathbb{P}\{T_a=k\} \, \mathbb{P}\{T_b=n-k\}.
\end{eqnarray*}
The same discussion holds also in the case $a \le -1$ and $b \le 0$.
\end{pf}

\begin{lem} The generating function of $T_a$ is given by
\begin{equation} \label{17:40}
\mathbb{E}\left(\xi^{T_a}\ind_{\{T_a <+\infty\}}\right)=\left\{\begin{array}{ll}\displaystyle{\left(\frac{1-\sqrt{1-4pq\xi^2}}{2q\xi}\right)^{\!a}} & \mbox{for $a \ge 1$,} \\[3ex] \displaystyle{\left(\frac{1-\sqrt{1-4pq\xi^2}}{2p\xi}\right)^{\!|a|}} & \mbox{for $a \le -1$.} \end{array}\right.
\end{equation}
\end{lem}
\begin{pf}
From Proposition \ref{Darling}, we have, for $a \ge 1$, $b \ge 0$ or $a \le -1$, $b \le 0$, that
$
\mathbb{P}\{S_n=a+b\}= \sum_{k=0}^{n} \mathbb{P}\{T_a=k\} \, \mathbb{P}\{S_{n-k}=b\}
$
which involves a convolution product. So, we derive
\begin{eqnarray*}
G(\xi, a+b)&=&\sum_{n=0}^{\infty} \mathbb{P}\{S_n=a+b\} \xi^n =\bigg[ \sum_{n=0}^{\infty} \mathbb{P}\{T_a=n\} \xi^n \bigg] \! \bigg[ \sum_{n=0}^{\infty} \mathbb{P}\{S_{n}=b\} \xi^n \bigg] \\
&=& \mathbb{E}\left(\xi^{T_a}\ind_{\{T_a <+\infty\}}\right) G(\xi, b)
\end{eqnarray*}
and then, for $b=0$,
$
\mathbb{E}\left(\xi^{T_a}\ind_{\{T_a <+\infty\}}\right)=\frac{G(\xi,a)}{G(\xi,0)},
$
which yields, thanks to (\ref{21:41}), result (\ref{17:40}).
\end{pf}

\begin{cor}
The distribution of $T_a$ is given as follows: for $1 \le |a| \le n$ such that $n+a$ is even,
\begin{equation} \label{18:04}
\mathbb{P}\{T_a=n\} =\frac{|a|}{n} \binom{n}{\frac{n+a}{2}} p^{\frac{n+a}{2}} q^{\frac{n-a}{2}}=\frac{|a|}{n} \, \mathbb{P}\{S_n=a\}
\end{equation}
and if $n+a$ is odd, this probability vanishes. Moreover, the probability of eventually hitting level $a$ is given by
\begin{equation} \label{18:04A}
\mathbb{P}\{T_a<+\infty\} =\left\{\begin{array}{ll}\displaystyle{\left(\frac{1-|p-q|}{2q}\right)^{\!a}} & \mbox{for $a \ge 1$,} \\[3ex] \displaystyle{\left(\frac{1-|p-q|}{2p}\right)^{\!|a|}} & \mbox{for $a \le -1$.} \end{array}\right.
\end{equation}
\end{cor}
\begin{pf}
From (\ref{17:42}), (\ref{21:41}) and (\ref{17:40}), we note that for $a \ge 1$,
\begin{eqnarray*}
\frac{\mathrm{d}}{\mathrm{d}\xi} \left[\mathbb{E}\left(\xi^{T_a}\ind_{\{T_a <+\infty\}}\right)\right] &=& \frac{a}{\xi \sqrt{1-4pq\xi^2}} \left(\frac{1-\sqrt{1-4pq\xi^2}}{2q\xi}\right)^{\!a}\\
&=&\frac{a}{\xi} \, G(\xi,a) = a \sum_{n=a}^{\infty} \binom{n}{\frac{n+a}{2}} p^{\frac{n+a}{2}} q^{\frac{n-a}{2}} \xi^{n-1}.
\end{eqnarray*}
By integrating this equality with respect to $\xi$, we obtain
$$
\mathbb{E}\left(\xi^{T_a} \ind_{\{T_a <+\infty\}}\right)=a \sum_{n=a}^{\infty} \binom{n}{\frac{n+a}{2}} p^{\frac{n+a}{2}} q^{\frac{n-a}{2}} \frac{\xi^{n}}{n}
$$
from which we extract (\ref{18:04}). The case $a \le -1$ can be treated exactly in the same way. Finally, (\ref{18:04A}) comes from (\ref{17:40}) by choosing $\xi =1$ therein and remarking that $1-4pq=(p-q)^2$.
\end{pf}

Notice the following equality which will be further used:
\begin{equation} \label{qq}
\mathbb{P} \{T_{-a}=n\}=\left(\frac{q}{p}\right)^{\!a} \mathbb{P} \{T_a=n\}.
\end{equation}

\begin{cor}
The distribution of $T_0$ is given as follows: for even $n \ge 1$,
$$
\mathbb{P}\{T_0=n\}= \binom{n}{\frac{n}{2}} \frac{(pq)^\frac{n}{2}}{n-1}
$$
and this probability vanishes for odd values of $n$. Moreover, the probability of eventually returning to the origin is given by
$
\mathbb{P}\{T_0<+\infty\}=1-|p-q|.
$
\end{cor}
\begin{pf}
By using the Markov property of the random walk as well as the translation invariance, we easily see that
$$
\mathbb{P}\{T_0=n\} = \mathbb{P}\{S_1=1\} \, \mathbb{P}\{T_{-1}=n-1\} + \mathbb{P}\{S_1=-1\} \, \mathbb{P}\{T_1=n-1\}.
$$
From (\ref{qq}), we observe that $p \, \mathbb{P}\{T_{-1}=n-1\}=q \, \mathbb{P} \{T_1=n-1\}$ and then
$$
\mathbb{P}\{T_0=n\} = 2q \, \mathbb{P} \{T_1=n-1\} = 2\binom{n-1}{\frac{n}{2}} \frac{(pq)^{\frac{n}{2}}}{n-1}= \binom{n}{\frac{n}{2}} \frac{(pq)^{\frac{n}{2}}}{n-1}.
$$
Thanks to the previous calculations and (\ref{17:40}) applied to $a=1$, we can write out the generating function of $T_0$:
\begin{eqnarray*}
\mathbb{E}\left(\xi^{T_0}\ind_{\{T_0 <+\infty\}}\right) &=& \sum_{n=1}^{\infty} \mathbb{P}\{T_0=n\} \xi^n = 2q \sum_{n=1}^{\infty} \mathbb{P}\{T_1=n-1\} \xi^n \\&=& 2q \xi \, \mathbb{E}\left(\xi^{T_1}\ind_{\{T_1 <+\infty\}}\right) = 1-\sqrt{1-4pq\xi^2}.
\end{eqnarray*}
For $\xi=1$, we immediately find the value of $\mathbb{P}\{T_0<+\infty\}$.
\end{pf}

\begin{rem}
In the case of the symmetric random walk ($p=q=\frac12$), the results of this subsection take the remarkable following form:
\begin{eqnarray*}
&\displaystyle\mathbb{P}\{S_n=j\} = \frac{1}{2^n} \binom{n}{\frac{n+j}{2}}, \qquad G(\xi,j) = \frac{1}{\sqrt{1-\xi^2}} \left(\frac{1- \sqrt{1-\xi^2}}{\xi}\right)^{\!|j|}\!, \\
&\hspace{-0.8cm}\!\!\displaystyle\mathbb{P}\{T_a=n\} = \frac{|a|}{n2^n} \binom{n}{\frac{n+a}{2}}\; \mbox{if } a\in\mathbb{Z}\setminus\{0\}, \; \mathbb{P}\{T_0=n\} = \frac{1}{(n-1)2^n} \binom{n}{\frac{n}{2}}, \; \mathbb{P}\{T_a<+\infty\} = 1\; \mbox{if }a\in\mathbb{Z}, \\
&  \displaystyle\mathbb{E}\left(\xi^{T_a}\right) = \left(\frac{1-\sqrt{1-\xi^2}}{\xi}\right)^{\!|a|}  a\in\mathbb{Z}\setminus\{0\}, \qquad \mathbb{E}\left(\xi^{T_0}\right) = 1-\sqrt{1-\xi^2}.
\end{eqnarray*}
\end{rem}

%%%%%%%%%%%%%%%%%%%%%%%%%%%%%
\subsection{Probabilistic interpretation of (\ref{s}) and (\ref{d})}

By rewriting the distributions of $S_n$ and $T_a$ (formulas (\ref{17:42}) and (\ref{18:04})) in terms of the Gamma function we have that
\begin{eqnarray}
\mathbb{P}\{S_n=i\}&=&\frac{\Gamma(n+1)}{\Gamma \left(\frac{n+i}{2}+1\right) \Gamma \left(\frac{n-i}{2}+1\right)} \, p^{\frac{n+i}{2}} q^{\frac{n-i}{2}}, \label{st} \\
\mathbb{P}\{T_a=j\}&=& |a| \, \frac{\Gamma(j)}{\Gamma \left(\frac{j+a}{2}+1\right) \Gamma \left(\frac{j-a}{2}+1\right)} \, p^{\frac{j+a}{2}} q^{\frac{j-a}{2}}. \label{st2}
\end{eqnarray}
Let $\mu$ and $\nu$ be positive integers. From (\ref{st}), with $i=\nu-1$ and $n=2r-2k+\nu-1$, we obtain
$$
\frac{\Gamma(2r-2k+\nu)}{(r-k)! \, \Gamma(r-k+\nu)}=p^{-(r-k+\nu-1)} q^{-(r-k)} \, \mathbb{P} \{S_{2r-2k+\nu-1}=\nu-1\},
$$
and, for $i=\mu+\nu-1$ and $n=2r+\mu+\nu-1$,
$$
\frac{\Gamma(2r+\mu+\nu)}{r! \,\Gamma(r+\mu+\nu)}=p^{-(r+\mu+\nu-1)} q^{-r} \, \mathbb{P}\{S_{2r+\mu+\nu-1}=\mu+\nu-1\}.
$$
From (\ref{st2}), with $a=\mu$ and $j=2k+\mu$, we also have
$$
\frac{\Gamma(2k+\mu)}{k! \, \Gamma(k+\mu+1)}=\frac{p^{-(k+\mu)} q^{-k} }{\mu} \,\mathbb{P}\{T_\mu=2k+\mu\}.
$$
Therefore, formula (\ref{s}) can be rewritten as a Darling-Siegert-type relationship (which extends formula (4.7) of \cite{D-S}):
$$
\sum_{k=0}^{r} \mathbb{P}\{T_\mu=2k+\mu\} \, \mathbb{P}\{S_{2r-2k+\nu-1}=\nu-1\}=\mathbb{P}\{S_{2r+\mu+\nu-1}=\mu+\nu-1\}.
$$
Similarly, by replacing successively in (\ref{st2}) $(j, a)$ by $(2k+\mu, \mu)$, $(2r-2k+\nu, \nu)$ and $(2r+\mu+\nu, \mu+\nu)$ we obtain the following expressions
\begin{eqnarray*}
\mathbb{P}\{T_\mu=2k+\mu\}&=&\mu\, \frac{\Gamma(2k+\mu)}{k! \, \Gamma( k+\mu+1)} \, p^{k+\mu} q^k, \\
\mathbb{P}\{T_\nu=2r-2k+\nu\}&=&\nu\, \frac{\Gamma(2r-2k+\nu)}{(r-k)! \, \Gamma(r- k+\nu+1)} \, p^{r-k+\nu} q^{r-k}, \\
\mathbb{P}\{T_{\mu+\nu}=2r+\mu+\nu\}&=& (\mu+\nu)\, \frac{\Gamma(2r+\mu+\nu)}{r! \, \Gamma(r+\mu+\nu+1)} \, p^{r+\mu+\nu} q^{r}.
\end{eqnarray*}
Formula (\ref{d}) now reads
$$
\sum_{k=0}^{r} \mathbb{P}\{T_\mu=2k+\mu\} \, \mathbb{P}\{T_\nu=2r-2k+\nu\} =\mathbb{P}\{T_{\mu+\nu}=2r+\mu+\nu\},
$$
which provides a probabilistic interpretation of (\ref{d}).\\

%%%%%%%%%%%%%%%%%%%%%%%%%%%%%
\subsection{Conditional aspect of (\ref{s}) and  (\ref{d})}

Let us divide equation (\ref{s}) by its right-hand side. This yields
$$
\sum_{k=0}^{r} \frac{r!}{k! \, (r-k)!} \, \frac{\Gamma(2k+\mu) \, \Gamma(2r-2k+\nu)}{\Gamma(2r+\mu+\nu)} \, \frac{\Gamma(r+\mu+\nu)}{\Gamma(k+\mu+1) \, \Gamma(r-k+\nu)} \mu=1,
$$
that is to say, by means of the binomial coefficients and changing $\nu$ into $\nu+1$, that for any integers $\mu \ge 1$ and $\nu \ge 0$,
$$
\sum_{k=0}^{r} \frac{\binom{r}{k} \binom{r+\mu+\nu}{k+\mu}}{\binom{2r+\mu+\nu}{2k+\mu}} \, \frac{\mu}{2k+\mu}=1.
$$
Then, we can observe that the family of numbers $(p_k^{\mu,\nu,r})_{0 \le k \le r}$ defined by
$
p_k^{\mu,\nu,r}= \frac{\binom{r}{k} \binom{r+\mu+\nu}{k+\mu}}{\binom{2r+\mu+\nu}{2k+\mu}} \, \frac{\mu}{2k+\mu}
$
is a probability distribution. Thanks to (\ref{17:42}) and (\ref{18:04}), we recognize that
%\begin{eqnarray*}
%p_k^{\mu,\nu,r}&=& \frac{r! \, (r+\mu+\nu)! \, (2k+\mu)! \, (2r-2k+\nu)!}{k! \, (r-k)! \, (k+\mu)! \, (r-k+\nu)! \, (2r+\mu+\nu)!} \, \frac{\mu}{2k+\mu} \\
%&=& \frac{r! \, (r+\mu+\nu)!}{(2r+\mu+\nu)!} \frac{(2k+\mu)!}{k! \, (k+\mu)!} \frac{(2r-2k+\nu)!}{(r-k)! \, (r-k+\nu)!} \, \frac{\mu}{2k+\mu} = \frac{\binom{2k+\mu}{k+\mu} \binom{2r-2k+\nu}{r-k+\nu}}{\binom{2r+\mu+\nu}{r+\mu+\nu}} \, \frac{\mu}{2k+\mu},
%\end{eqnarray*}
%we recognize, thanks to (\ref{17:42}) and (\ref{18:04}), that
$$
p_k^{\mu,\nu,r}=\frac{\mathbb{P} \{T_\mu=2k+\mu\} \, \mathbb{P} \{S_{2r-2k+\nu}=\nu\}}{\mathbb{P} \{S_{2r+\mu+\nu}=\mu+\nu\}}.
$$
Next, by the Markov property and the translation invariance, we have
\begin{eqnarray*}
 \lefteqn{\mathbb{P} \{T_\mu=2k+\mu\} \, \mathbb{P}\{S_{2r-2k+\nu}=\nu\}=\mathbb{P} \{T_\mu=2k+\mu\} \, \mathbb{P}\{S_{2r+\mu+\nu}-S_{2k+\mu}=\nu\} }\\
&=& \mathbb{P} \{T_\mu=2k+\mu, S_{2r+\mu+\nu}-S_{2k+\mu}=\nu\} =\mathbb{P} \{T_\mu=2k+\mu,S_{2r+\mu+\nu}=\mu+\nu\}
\end{eqnarray*}
which leads to the following probabilistic interpretation for $p_k^{\mu,\nu,r}$:
$$
p_k^{\mu,\nu,r}= \mathbb{P} \{T_\mu=2k+\mu | S_{2r+\mu+\nu}=\mu+\nu\}.
$$
This means that $(p_k^{\mu,\nu,r})_{0 \le k \le r}$ is the distribution of the first hitting time of level $\mu$ for the random walk pinned at the extremity $\mu+\nu$ at time $2r+\mu+\nu$. Of course, this interpretation can also be directly derived from (\ref{d}).

Let us consider the particular case $\mu = 1$ and $\nu = 0$ and set $q^r_k=p_{k-1}^{1,0,r-1}$ for $1 \le k \le r$. We have
\begin{equation}\label{o}
q_k^r= \mathbb{P} \{T_1=2k-1| S_{2r-1}=1\}.
\end{equation}
We also have
\begin{equation}\label{oo}
q_k^r= \mathbb{P} \{T_{-1}=2k-1| S_{2r-1}=-1\}.
\end{equation}
Indeed, by (\ref{qqq}) and (\ref{qq}),
\begin{eqnarray*}
\mathbb{P} \{T_{-1}=2k-1 | S_{2r-1}=-1\}&=&\frac{\mathbb{P} \{T_{-1}=2k-1\} \, \mathbb{P}\{S_{2r-2k}=0\}}{\mathbb{P}\{S_{2r-1}=-1\}}\\
&=& \frac{\frac{q}{p} \, \mathbb{P} \{T_1=2k-1\} \, \mathbb{P}\{S_{2r-2k}=0\}}{\frac{q}{p} \, \mathbb{P}\{S_{2r-1}=1\}}\\
&=&\mathbb{P} \{T_1=2k-1| S_{2r-1}=1\} = q_k^r.
\end{eqnarray*}
Actually, our aim is to check that
\begin{equation}\label{ox}
q_k^r=\mathbb{P}\{T_0=2k|S_{2r}=0\}.
\end{equation}
In this form, we can see that the probability distribution $(q_k^r)_{1 \le k \le r}$ is nothing but that of the first return time to $0$ for the bridge of length $2r$ of the random walk. This distribution appeared in \cite{D-S} in the symmetric case $p=q=\frac{1}{2}$ and was a consequence of (\ref{d}) whereas, here, we use (\ref{s}) instead. Let us now verify (\ref{ox}). For this, we write
\begin{eqnarray*}
\mathbb{P} \{T_0=2k, S_{2r}=0\}&=&\mathbb{P}\{T_0=2k\} \, \mathbb{P}\{S_{2r-2k}=0\} \\
&=& \left(p \, \mathbb{P} \{T_{-1}=2k-1\}+q \, \mathbb{P} \{T_1=2k-1\}\right) \mathbb{P}\{S_{2r-2k}=0\} \\
&=&p \, \mathbb{P}\{T_{-1}=2k-1, S_{2r-1}=-1\} + q \, \mathbb{P}\{T_1=2k-1, S_{2r-1}=1\} \\
&=&p \, \mathbb{P}\{T_{-1}=2k-1| S_{2r-1}=-1\} \, \mathbb{P}\{S_{2r-1}=-1\} \\
&&+\, q \, \mathbb{P} \{T_1=2k-1| S_{2r-1}=1\} \, \mathbb{P}\{S_{2r-1}=1\}.
\end{eqnarray*}
By using (\ref{o}), (\ref{oo}), we have $\mathbb{P}\{T_0=2k, S_{2r}=0\}= q_k^r \left(p \, \mathbb{P} \{S_{2r-1}=-1\}+q \, \mathbb{P}\{S_{2r-1}=1\}\right) =q_k^r \,\mathbb{P}\{S_{2r}=0\}$, which immediately proves (\ref{ox}).

%%%%%%%%%%%%%%%%%%%%%%%%%%%%%
\subsection{Probabilistic interpretation of (\ref{b})}\label{subsection}

In the spirit of this paper, we give a probabilistic interpretation for the simple formula (\ref{b}). For this, in view of (\ref{serie}), we have the equality $\sum_{n=0}^{\infty} \frac{\Gamma(n+\mu)}{n! \, \Gamma(\mu)}\, p^n q^\mu=1$ where $p \in (0,1)$ is fixed and $q=1-p$. This suggests to introduce a family of random variables $(\overline{T}_\mu)_{\mu \in \mathbb{N} \setminus \{0\}}$, with probability distribution
$$
\mathbb{P}\{\overline{T}_\mu=n\}=\frac{\Gamma(n+\mu)}{n! \, \Gamma(\mu)} \, p^n q^\mu, \qquad n \in \mathbb{N}.
$$
This is the well-known negative binomial distribution with parameters $(\mu, p)$. Formula (\ref{b}) implies that this family satisfies, for any $r \in \mathbb{N}$, the relationship
$$
\sum_{k=0}^{r} \mathbb{P} \{\overline{T}_\mu=k\} \, \mathbb{P}\{\overline{T}_\nu=r-k\}=\mathbb{P}\{\overline{T}_{\mu+\nu}=r\}.
$$
This formula can be rephrased by introducing an independent copy $(\overline{T}'_\mu)_{\mu>0}$ of the family $(\overline{T}_\mu)_{\mu>0}$ as $\mathbb{P} \{\overline{T}_\mu+ \overline{T}'_\nu=r\}=\mathbb{P}\{\overline{T}_{\mu+\nu}=r\}$, for any $r \in \mathbb{N}$. In words, for any $\mu, \nu >0$, the random variables $\overline{T}_\mu+\overline{T}'_\nu$ and $\overline{T}_{\mu+\nu}$ have the same probability distribution. If we choose the random variables $(\overline{T}_\mu)_{\mu>0}$ independent, then, for any different $\mu, \nu>0$, $\overline{T}_\mu+\overline{T}_\nu$ and $\overline{T}_{\mu+\nu}$ have the same laws.

Another interpretation consists of relating the family $(\overline{T}_{\mu})_{\mu>0}$ to the family of passage times for a certain random walk $(S'_n)_{n \in \mathbb{N}}$, which we are now going to describe. We can rewrite (\ref{b}) in the following form: for positive integers $\mu, \nu$,
\begin{equation}\label{binom}
\sum_{k=0}^{r} \binom{k+\mu-1}{\mu-1} \binom{r-k+\nu-1}{\nu-1}=\binom{r+\mu+ \nu-1}{\mu+\nu-1}.
\end{equation}
This formula suggests to introduce the following family of random variables. Fix $p \in (0,1)$, $q=1-p$ and let $\left(T'_\mu\right)_{\mu \in \mathbb{N} \setminus \{0\} }$ be the Pascal random variable with parameters $(\mu,p)$ (this is the negative binomial distribution translated by $\mu$):
\begin{equation}\label{binomdist}
\mathbb{P}\{T'_\mu=n\}=\binom{n-1}{\mu-1} p^{\mu} q^{n-\mu}, \hspace{1cm} n \ge \mu.
\end{equation}
Formula (\ref{binom}) reads
$$
\sum_{k=0}^{r} \mathbb{P}\{T'_\mu=k+\mu\} \, \mathbb{P}\{T'_\nu=r-k+\nu\} = \mathbb{P}\{T'_{\mu+\nu}=r+\mu+\nu\}.
$$

Our goal is to find a certain process $\left(S'_n\right)_{n \in \mathbb{N}}$ living on $\mathbb{N}$ that obeys the strong Markov property, that is translation invariant, and such that for any positive integer $\mu$ we have $T'_\mu=\min\{n \ge 1: S'_n = \mu\}$ (in order to have $S'_{T'_\mu}=\mu$). In view of Proposition \ref{Darling}, for any $a \in \mathbb{N}$, the generating function $G'(\xi,a)$ defined, for $|\xi|<1$, by $G'(\xi,a)=\sum_{n=0}^{\infty} \mathbb{P}\{S'_n=a\} \xi^n$ satisfies the following relation
\begin{equation}\label{generat}
G'(\xi,a)=G'(\xi,0)\,\mathbb{E}(\xi^{T'_a}).
\end{equation}
Above, the generating function of $T'_a$ is explicitly given, with the aid of (\ref{serie}), by
\begin{equation}\label{generatbis}
\mathbb{E}(\xi^{T'_a})=\sum_{n=a}^{\infty} \binom{n-1}{a-1} p^a q^{n-a} \xi^n = \sum_{n=0}^{\infty} \binom{n+a-1}{a-1} p^a q^{n} \xi^{n+a} =\left(\frac{p\,\xi}{1-q\,\xi}\right)^{\!a}
\end{equation}
and we have
\begin{equation}\label{generatter}
G'(\xi,0)=\frac{1}{1-q \xi}
\end{equation}
since, by observing that we obviously must have $\sum_{a=0}^{\infty} \mathbb{P}\{S'_n=a\}=1$ and by referring again to (\ref{serie}),
\begin{eqnarray*}
\hspace{-0.7cm}\frac{1}{1-\xi}&=& \sum_{n=0}^{\infty} \xi^n =\sum_{n=0}^{\infty} \sum_{a=0}^{\infty} \mathbb{P}\{S'_n=a\} \xi^n = \sum_{a=0}^{\infty} G' (\xi,a) =G'(\xi,0) \sum_{a=0}^{\infty} \left(\frac{p\, \xi}{1-q\,\xi}\right)^{\!a} = G'(\xi,0) \, \frac{1-q\xi}{1-\xi}.
\end{eqnarray*}
By plugging (\ref{generatbis}) and (\ref{generatter}) into (\ref{generat}) and using (\ref{serie}), we obtain that
$$
G'(\xi,a)=\frac{(p \xi)^a}{(1-q\xi)^{a+1}} = (p \xi)^a \sum_{n=0}^{\infty} \binom{n+a}{a} q^n \xi^n = (p \xi)^a \sum_{n=a}^{\infty} \binom{n}{a} q^{n-a} \xi^{n-a} = \sum_{n=a}^{\infty} \binom{n}{a} p^a q^{n-a} \xi^n
$$
which implies that
\begin{equation}\label{binomdistbis}
\mathbb{P}\{S'_n=a\}=\binom{n}{a} p^{a} q^{n-a}, \hspace{2cm} 0 \le a \le n.
\end{equation}
So, we recognize the binomial distribution of parameters $(n,p)$. We deduce from this discussion that the process we were seeking is the non-decreasing random walk $(S'_n)_{n \in \mathbb{N}}$ where $S'_n=\sum_{k=0}^{n} X'_k$ with $X'_0=0$ and such that the $X'_k$, $k \ge 1$, are independent random variables taking the value 1 or 0 with probability $p$ or $q$, respectively.

Conversely, we can easily check that this random walk satisfies the three assumptions (strong Markov property, translation invariance, discrete continuity) made in the above analysis.

This random walk plays an important role in game theory. Indeed, if we toss $n$ times a coin and call ``head" a success, $S'_n$ is the number of successes among the $n$ tosses while $T'_a$ is the minimal number of tosses for obtaining exactly $a$ successes.

We have interpreted (\ref{b}) by means only of the variables $(T'_a)_{a \in \mathbb{N} \setminus \{0\}}$. Since we have identified the underlying random walk $(S'_n)_{n\in \mathbb{N}}$ associated with $(T'_a)_{a \in \mathbb{N} \setminus \{0\}}$, let us focus ourselves on the Darling-Siegert relationship connecting both families. We have $\mathbb{P}\{S'_n=a+b\}=\sum_{k=a}^{n-b} \mathbb{P}\{T'_a=k\} \, \mathbb{P}\{S'_{n-k}=b\}$
which, by (\ref{binomdist}) and (\ref{binomdistbis}) supplies the simple combinatorial identity
$$
\sum_{k=a}^{n-b} \binom{k-1}{a-1} \binom{n-k}{b}=\binom{n}{a+b}.
$$
Rewriting this latter, by setting $n=r+a+b$, as
$$
\sum_{k=0}^{r} \binom{k+a-1}{a-1} \binom{r-k+b}{b}=\binom{r+a+b}{a+b},
$$
we retrieve, by translating one of the parameters $a$ and $b$, formula (\ref{binom}). We then conclude that both Darling-Siegert-type relationships related to $(S'_n)_{n\in \mathbb{N}}$ together with $(T'_a)_{a \in \mathbb{N} \setminus \{0\}} $
$$\mathbb{P}\{S'_n=a+b\}=\sum_{k=a}^{n-b} \mathbb{P}\{T'_a=k\} \, \mathbb{P}\{S'_{n-k}=b\}, \hspace{0.5cm} \mathbb{P}\{T'_{a+b}=n\}=\sum_{k=a}^{n-b} \mathbb{P}\{T'_a=k\} \, \mathbb{P}\{T'_b=n-k\} $$
lead, in this case, to the same combinatorial identity, contrarily to the foregoing case. This fact is in good accordance with the observation that both identities (\ref{16:40A}) and (\ref{16:40B}) lead to the same combinatorial identity (\ref{b}).

\begin{rem}\label{remark}
The particular relationship $ \mathbb{P}\{T_a=n\}=\frac{a}{n} \, \mathbb{P}\{S_n=a\}$, for any $a \in \mathbb{N} \setminus \{0\}$, which holds for both random walks considered in this paper, shows that one Darling-Siegert-type formula implies the other one. Indeed, starting from
$
\sum_{k=0}^{n} \mathbb{P}\{T_a=k\} \, \mathbb{P}\{S_{n-k}=b\}= \mathbb{P}\{S_n=a+b\},
$
we obtain
$
\sum_{k=0}^{n} (n-k) \, \mathbb{P}\{T_a=k\} \, \mathbb{P}\{T_b=n-k\}=\frac{bn}{a+b} \, \mathbb{P}\{T_{a+b}=n\}.
$
By making the change of index $k\to n-k$ and by interchanging the roles of $a$ and $b$, this can be written as
$
\sum_{k=0}^{n} k \, \mathbb{P}\{T_a=k\} \, \mathbb{P}\{T_b=n-k\}=\frac{an}{a+b} \, \mathbb{P}\{T_{a+b}=n\}.
$
Now, adding the two foregoing equalities leads to formula
$
\sum_{k=0}^{n} \mathbb{P}\{T_a=k\} \, \mathbb{P}\{T_b=n-k\}= \mathbb{P}\{T_{a+b}=n\}.
$
\end{rem}

% \appendix

% \section{}
% \label{}

% Bibliographic references with the natbib package:
% Parenthetical: \citep{Bai92} produces (Bailyn 1992).
% Textual: \citet{Bai95} produces Bailyn et al. (1995).
% An affix and part of a reference:
%   \citep[e.g.][Ch. 2]{Bar76}
%   produces (e.g. Barnes et al. 1976, Ch. 2).

\end{document}